\documentclass[12 pt]{amsart}
\usepackage{amscd, amssymb}
\newtheorem{tm}{Theorem}
\newtheorem{co}{Conjecture}

\newcommand{\Q}{\mathbb{Q}}

\newcommand{\la}{\lambda}

\newcommand{\mbar}[1]{{\overline{M}}_{#1}}
\newcommand{\mct}[1]{M^{ct}_{#1}}
\newcommand{\mrt}[1]{M^{rt}_{#1}}

\newcommand{\rbar}[1]{R^{\bullet}({\overline{M}}_{#1})}
\newcommand{\rct}[1]{R^{\bullet}(M^{ct}_{#1})}
\newcommand{\rrt}[1]{R^{\bullet}(M^{rt}_{#1})}
\newcommand{\rsm}[1]{R^{\bullet}(M_{#1})}

\newcommand{\rcct}[1]{R^{\bullet}_c(M^{ct}_{#1})}
\newcommand{\rcrt}[1]{R^{\bullet}_c(M^{rt}_{#1})}
\newcommand{\rcsm}[1]{R^{\bullet}_c(M_{#1})}

\newcommand{\bpf}{\noindent {\em Proof.} }
\newcommand{\epf}{\qed \vspace{+10pt}}

\begin{document}
\title{A remark on a conjecture of\\ Hain and Looijenga}
\author{Carel Faber}
\address{Department of Mathematics, KTH Royal Institute of Technology,
Lindstedtsv\"agen 25, 10044 Stockholm, Sweden.}
\email{faber@math.kth.se}
\date{}
\maketitle

\pagestyle{plain}
Let~$M_{g,n}$ (resp.~$\mbar{g,n}$) be the moduli space of 
smooth (resp.~stable) $n$-pointed curves of
genus~$g$ and let~$\mct{g,n}$ be the moduli space of pointed curves of
compact type, the complement of the boundary divisor~$\Delta_{\text{irr}}$
of irreducible singular curves and their degenerations. Let~$\mrt{g,n}$
be the moduli space of pointed curves with rational tails; for~$g\ge2$,
it is the inverse image of~$M_g$ under the natural 
morphism~$\mbar{g,n}\to\mbar{g}$, while~$\mrt{1,n}=\mct{1,n}$
and~$\mrt{0,n}=\mbar{0,n}$ by definition. Here, $(g,n)$ is a pair
of nonnegative integers such that~$2g-2+n>0$. There is a natural partial
ordering of these pairs: $(h,m)\le(g,n)$ if and only if~$h\le g$
and~$2h-2+m\le 2g-2+n$, or, in other words, if and only if there exists
a stable~$n$-pointed curve of genus~$g$ whose dual graph contains
a vertex of genus~$h$ with valency~$m$.

We recall the definition of the tautological algebras~$\rbar{g,n}$
from [FP2]: the system~$\{\rbar{g,n}\}_{(g,n)}$ is defined as the
set of smallest~$\Q$-subalgebras of the rational Chow
rings~$A^{\bullet}(\mbar{g,n})$ that is closed under push-forward
via all maps forgetting markings and all standard gluing maps.
The well-known~$\psi$-, $\kappa$-, and~$\lambda$-classes are
tautological. The system is also closed under pull-back via the
forgetting and gluing maps. The successive quotients~$\rct{g,n}$, $\rrt{g,n}$,
and~$\rsm{g,n}$ are defined as the restrictions to the respective
open subsets. (Observe that it is in general not known whether
the corresponding tautological localization sequences are exact in
the middle.)

The following results are known:
\begin{enumerate}
\item[(a)] $\rrt{g,n}$ vanishes in degrees~$>g-2+n-\delta_{0g}$
and is~$1$-dimensional in degree~$g-2+n-\delta_{0g}$.
\item[(b)] $\rct{g,n}$ vanishes in degrees~$>2g-3+n$
and is~$1$-dimensional in degree~$2g-3+n$.
\item[(c)] $\rbar{g,n}$ (vanishes in degrees~$>3g-3+n$
and) is~$1$-dimensional in degree~$3g-3+n$.
\end{enumerate}
Statement (a) was proved by Looijenga [L] and Faber [F1], [FP1]. 
Statements (b) and (c)
were proved by Graber and Vakil [GV1], [GV2] and Faber and Pandharipande
[FP2].

Recall the following three conjectures:
\begin{enumerate}
\item[(A)] $\rrt{g,n}$ is Gorenstein with socle in degree~$g-2+n-\delta_{0g}$.
\item[(B)] $\rct{g,n}$ is Gorenstein with socle in degree~$2g-3+n$.
\item[(C)] $\rbar{g,n}$ is Gorenstein with socle in degree~$3g-3+n$.
\end{enumerate}
(For a graded $\Q$-algebra $R^{\bullet}$, to be Gorenstein with socle
in degree~$m$ means that it vanishes in degrees $>m$, that $R^m$ is
isomorphic to $\Q$, and that the pairings $R^i\times R^{m-i}\to R^m$
are perfect.)

In the case~$g=0$, the three conjectures coincide and have been proved by Keel [K].
Conjecture (A) in the case~$n=0$ is due to the author [F1] and is true
for~$g\le23$. Hain and Looijenga [HL] raised (C) as a question
and (A), (B), and (C) were formulated in [Pa] (see also [FP1], [F2]).

Hain and Looijenga also introduce a compactly supported version
of the tautological algebra: they define~$\rcsm{g,n}$ as the set of
elements in~$\rbar{g,n}$ that restrict trivially to the Deligne-Mumford
boundary (i.e., the pull-back via any standard map from a product
of moduli spaces~$\mbar{g_i,n_i}$ onto the closure of a boundary
stratum vanishes). It is a graded ideal in~$\rbar{g,n}$ and a module
over~$\rsm{g,n}$. They then formulate the following conjecture
in the case~$n=0$:
\begin{co}[Hain and Looijenga~\cite{HL}] \label{HL}
The intersection pairings
$$
R^k(M_g)\times R^{3g-3-k}_c(M_g)\to R^{3g-3}_c(M_g)\cong\Q ,
\qquad k=0,1,2,\dots
$$
are perfect (Poincar\'e duality) and~$\rcsm{g}$
is a free~$\rsm{g}$-module of rank one.
\end{co}
Observe that~$\la_1\in R^1_c(M_{1,1})$ and~$\la_g\la_{g-1}\in R^{2g-1}_c(M_g)$
for~$g>1$. (The author's proof of the nonvanishing of~$R^{g-2+n}(\mrt{g,n})$
for~$g>0$ uses this fact.) So this
class is supposed to be a generator of the~$\rsm{g}$-module~$\rcsm{g}$
(the unique generator of degree~$2g-1$ up to a scalar).

However, the pull-backs of these classes to~$\mbar{g,n}$ don't lie
in~$\rcsm{g,n}$ for~$n\ge2$, since they don't vanish on the boundary
strata corresponding to curves with rational tails. Let us therefore
define~$\rcrt{g,n}$ as the set of elements in~$\rbar{g,n}$ that restrict 
trivially to~$\mbar{g,n}\setminus\mrt{g,n}$. Consider the following
conjectures:
\begin{enumerate}
\item[(D)] The intersection pairings
$$
R^k(\mrt{g,n})\times R^{3g-3+n-k}_c(\mrt{g,n})\to R^{3g-3+n}_c(\mrt{g,n})\cong\Q
$$
are perfect for $k\ge0$.
\item[(E)] In addition to (D), $\rcrt{g,n}$ is a free~$\rrt{g,n}$-module
of rank one.
\end{enumerate}
Conjecture (E) appears to be the natural generalization
of Conjecture~\ref{HL} to the case~$n>0$. For reasons that will become clear
in a moment, we also include the weaker statement~(D).
Observe that (E) implies that~$\la_g\la_{g-1}$ is a generator of~$\rcrt{g,n}$
for~$g>0$ (the unique one of degree~$2g-1$ up to a scalar), by (a) above.
\begin{tm} Conjectures (A) and (C) are true for all~$(g,n)$ if and only if
Conjecture (E) is true for all~$(g,n)$. More precisely,
$$A_{(g,n)}\,\,\,{\rm and}\,\,\, C_{(g,n)}
\,\Rightarrow\, E_{(g,n)}
\,\Rightarrow\, A_{(g,n)}\,\,\,{\rm and}\,\,\, D_{(g,n)}
$$
and
$$\{ D_{(g',n')} \}_{(g',n')\le(g,n)} \,\Rightarrow\,
\{ C_{(g',n')} \}_{(g',n')\le(g,n)}\,.$$
\end{tm}

\bpf Suppose first that (C) is not true for all~$(g,n)$ and let
a minimal counterexample be given by~$0\neq\alpha\in\rbar{g,n}$, i.e.,
$\rbar{g',n'}$ is Gorenstein for all~$(g',n')<(g,n)$ 
and $\deg(\alpha\beta)=0$ for all~$\beta\in\rbar{g,n}$. 
(We write $\deg$ for the degree homomorphism on~$R_0(\mbar{g,n})$
and its extension by zero to all of~$\rbar{g,n}$).
It follows that~$g>0$.

Let~$\pi$ denote the standard map~$\mbar{g-1,n+2}\to\mbar{g,n}$ 
onto the boundary divisor~$\Delta_{\text{irr}}$. 
Let~$\gamma\in\rbar{g-1,n+2}$ be arbitrary. Then 
$$\deg((\pi^*\alpha)\gamma)=\deg(\pi_*((\pi^*\alpha)\gamma))
=\deg(\alpha\pi_*\gamma)=0,$$
since $\pi_*\gamma$ is tautological.
Since~$\rbar{g-1,n+2}$ is Gorenstein, it follows that~$\pi^*\alpha=0$.

Next, let~$\pi$ denote one of the
standard maps~$\mbar{g_1,n_1}\times\mbar{g_2,n_2}\to\mbar{g,n}$ onto
a boundary component parametrizing reducible singular curves ($g_1+g_2=g$
and~$n_1+n_2=n+2$). We have the push-forward map
$$\pi_*:\rbar{g_1,n_1}\otimes_{\Q}\rbar{g_2,n_2}\to\rbar{g,n}$$
and the pull-back map in the other direction (cf.~[GP]).
The tensor product is Gorenstein, with perfect pairing given by
$$\deg((\beta_1\otimes\beta_2)(\gamma_1\otimes\gamma_2))
=\deg(\beta_1\gamma_1)\deg(\beta_2\gamma_2).$$
Let~$\gamma_1$ resp.~$\gamma_2$ be arbitrary elements
of~$\rbar{g_1,n_1}$ resp.~$\rbar{g_2,n_2}$. Then
$$\deg((\pi^*\alpha)(\gamma_1\otimes\gamma_2))
=\deg(\pi_*((\pi^*\alpha)(\gamma_1\otimes\gamma_2)))
=\deg(\alpha\pi_*(\gamma_1\otimes\gamma_2))=0,$$
since~$\pi_*(\gamma_1\otimes\gamma_2)$ is tautological. Again, it follows
that~$\pi^*\alpha=0$.

Therefore, $0\neq\alpha\in\rcsm{g,n}$ and a fortiori~$0\neq\alpha\in\rcrt{g,n}$.
But it pairs to zero with all~$\beta$ and this contradicts~$D_{(g,n)}$.
The implication in the second display follows as an immediate consequence.

The next step is to prove the implication~$E_{(g,n)}\Rightarrow A_{(g,n)}$.
As mentioned above, if~$g>0$ and~$E_{(g,n)}$ holds, then~$\la_g\la_{g-1}$ 
generates~$\rcrt{g,n}$ freely.
Suppose that~$A_{(g,n)}$ fails: let~$0\neq\alpha\in\rrt{g,n}$ be such that
it pairs to zero with all~$\beta\in\rrt{g,n}$, i.e., 
$\deg(\alpha\beta\la_g\la_{g-1})=0$ for all~$\beta$ (note that~$g>0$).
From~$D_{(g,n)}$, it follows that~$\alpha\la_g\la_{g-1}=0$, but this
contradicts~$E_{(g,n)}$. This proves the second implication in the first display.

To prove the first implication, we first show that~$A_{(g,n)}$
and~$C_{(g,n)}$ imply~$D_{(g,n)}$. Assume that~$D_{(g,n)}$ fails;
the perfect pairing may fail on either side.
Suppose first that~$0\neq\alpha\in\rcrt{g,n}$ pairs to zero with all
of~$\rrt{g,n}$. We know that~$\pi^*\alpha=0$, for every standard
map~$\pi$ associated to a stratum in~$\mbar{g,n}\setminus\mrt{g,n}$.
This means that the product of~$\alpha$ and a Chow class
pushed forward via such a map is zero (hence the pairing is well-defined).
Since~$\alpha$ pairs to zero with all of~$\rrt{g,n}$, it gives a counterexample
to~$C_{(g,n)}$. If instead~$0\neq\alpha\in\rrt{g,n}$ pairs to zero
with all of~$\rcrt{g,n}$, then it pairs to zero with all classes
of the form~$\beta\la_g\la_{g-1}$, for~$\beta\in\rrt{g,n}$ (note that~$g>0$).
In this case, $\alpha$ gives a counterexample to~$A_{(g,n)}$.

We conclude by showing that~$A_{(g,n)}$ and~$C_{(g,n)}$ imply~$E_{(g,n)}$.
We already have~$D_{(g,n)}$. If~$E_{(g,n)}$ doesn't hold, then~$g>0$
and certainly~$\la_g\la_{g-1}$ fails to be a basis for~$\rcrt{g,n}$, i.e.,
multiplication by~$\la_g\la_{g-1}$ fails to be surjective or injective.
From~$A_{(g,n)}$ and~$D_{(g,n)}$, it follows that the surjectivity and
injectivity of this map are equivalent (recall from [GP], Cor.~1,
that~$\rbar{g,n}$ is finite-dimensional). But if~$0\neq\alpha\in\rrt{g,n}$
and~$\alpha\la_g\la_{g-1}=0$, then~$A_{(g,n)}$ fails. \epf

There is an analogous result in the compact type case. Begin by defining~$\rcct{g,n}$
as the set of elements in~$\rbar{g,n}$ that pull back to zero via the
standard map~$\mbar{g-1,n+2}\to\mbar{g,n}$ onto~$\Delta_{\text{irr}}$.
Conjectures (D) and (E) have obvious analogues $(\text{D}^{\text{ct}})$
and $(\text{E}^{\text{ct}})$. We have that
$$B_{(g,n)}\,\,\,{\rm and}\,\,\, C_{(g,n)}
\,\Rightarrow\, E^{ct}_{(g,n)}
\,\Rightarrow\, B_{(g,n)}\,\,\,{\rm and}\,\,\, D^{ct}_{(g,n)}
$$
and
$$\{ D^{ct}_{(g',n')} \}_{(g',n')\le(g,n)} \,\Rightarrow\,
\{ C_{(g',n')} \}_{(g',n')\le(g,n)}\,.$$
The proof proceeds entirely analogously; the class~$\la_g$ now plays
the role of~$\la_g\la_{g-1}$ (it is no longer necessary to treat
the case~$g=0$ separately).

\medskip
{\bf Acknowledgements.\/} The author thanks Eduard Looijenga, Rahul
Pandharipande, and Michael Shapiro for useful discussions.
The author is supported by 
the G\"oran Gustafsson Foundation for Research in Natural Sciences
and Medicine and grant 622-2003-1123 from the Swedish Research Council.

\medskip
{\em Note added in the second version\/} (November 2010).
Tavakol [T1] has proved that the tautological ring of $\mct{1,n}
=\mrt{1,n}$ is Gorenstein with socle in degree $n-1$
(Conjectures (A) and (B) for $g=1$). From Theorem~$1$, the tautological rings
$\rbar{1,n}$ are Gorenstein if and only if $E_{(1,n)}$ holds for all~$n\ge1$,
in other words, if and only if 
$\rcrt{1,n}$ is 
generated by~$\la_1$
as an $\rrt{1,n}$-module.

\medskip
{\em Note added in the third version\/} (April 2012).
Tavakol [T2] has now also proved Conjecture (A) for $g=2$:
the tautological ring of 
$\mrt{2,n}$ is Gorenstein with socle in degree $n$.

\medskip
{\em Note added in the fourth version\/} (June 2012).
Petersen [Pe] has proved that the tautological ring of $\mbar{1,n}$
is Gorenstein with socle in degree $n$
(Conjecture (C) for $g=1$).

\end{document}